\documentclass[12pt,a4paper]{ijicic}
\usepackage{graphicx}
\usepackage{amsmath,amssymb}

\def\currentvolume{5}
\def\currentissue{x}
\def\currentyear{2009}
\def\currentmonth{x}
\def\ppages{0--0}

\makeatletter
\DeclareMathOperator{\erf}{erf}
\def\lfig#1{\label{fig:#1}}\def\rfig#1{Fig.\ref{fig:#1}}
\def\Rfig#1{Figure \ref{fig:#1}}
\def\leqn#1{\label{eqn:#1}}\def\reqn#1{(\ref{eqn:#1})}

\let\gv=\bi
\def\dgv#1{\dot{\gv{#1}}}
\def\@pd[#1]#2#3{{\def\aaa{#1}\ifx\aaa\@empty%
 \frac{\partial #2}{\partial {#3}}%
 \else\frac{\partial^{#1} #2}{\partial {#3}^{#1}}\fi}}
\def\@od[#1]#2#3{{\def\aaa{#1}\ifx\aaa\@empty%
 \frac{d #2}{d #3}%
 \else\frac{d^{#1} #2}{d {#3}^{#1}}\fi}}
\def\pd{\@ifnextchar[{\@pd}{\@pd[]}}
\def\od{\@ifnextchar[{\@od}{\@od[]}}
\def\ave#1{\langle#1\rangle}
\def\lam{\lambda}
\makeatother

\title{Bifurcation Analysis of Noise-induced Synchronization}

\begin{document}
\maketitle

\centerline{\scshape  Katsutoshi Yoshida}
 \medskip
{\footnotesize
\centerline{Department of Mechanical Systems Engineering}
\centerline{Utsunomiya University} 
\centerline{7-1-2 Yoto, Utsunomiya-shi, Tochigi 321-8585, Japan}
\centerline{yoshidak@cc.utsunomiya-u.ac.jp} }
\medskip

\centerline{\scshape  Yusuke Nishizawa}
 \medskip
{\footnotesize
\centerline{ShinMaywa Industries, Ltd.} 
\centerline{1-1 Shinmeiwa-cho, Takarazuka-shi, Hyogo 665-8550, Japan} }
\medskip

\begin{abstract}
We investigate bifurcation phenomena between slow and fast
convergences of synchronization errors arising in the proposed
synchronization system consisting of two identical nonlinear
dynamical systems linked by a common noisy input only.  The
numerical continuation of the saddle-node bifurcation set of the
primary resonance of moments provides an effective identifier of the
slow convergence of synchronization errors.
\\

\end{abstract}

\section{Introduction}

The noise-induced synchronization of a dynamical system with its 
copies can easily be found in nonlinear systems, such as the 
discrete maps\cite{toral01,suetani04}, the Lorenz
system\cite{toral01}, the Duffing oscillator\cite{stefan03}, the 
single mode CO$_2$ laser\cite{zhou03}, and the uncoupled
neurons\cite{neiman02}.  One of the most important results of them
is that the perfect synchronization may arise under some suitable
conditions\cite{toral01,stefan03,zhou03}.  Moreover, the perfect
synchronization exhibits significant degree of robustness against
mismatches among the copies such as the parameters
mismatch\cite{stefan03} and the independent random fluctuations of
the copies\cite{toral01,yoshida04}.  Furthermore, regarding the 
response of the synchronization system as a Markov process to derive the 
transition law of it, we have analytically shown that the perfect
synchronization can be regarded as an absorbing barrier of the 
Markov process\cite{yoshida06}. 

In these studies, however, little research have been done on
transient behavior to approach the perfect synchronization.  In
engineering applications, too slow convergence of synchronization
errors would be regarded as failing to converge, even if the
synchronization is achieved mathematically.
In practice, the perfect synchronization is possibly applicable to
synchronizing initial conditions of coexisting oscillators with a common
specification such as independent subcircuits in a circuit system,
independent mechanical vibrators on bench testing, and so on.  However,
it is hardly applicable to industrial purposes when the convergence
speed is too slow.

To solve this problem, we have already investigated how to characterize
the slow convergence of synchronization
errors\cite{sss06,yoshida07}, showing that the slow convergence is
related to the slow diffusions caused by multimodal probability
densities so that it can be detected as multi-valued solutions of
the moment differential equations (MDE).

In this paper, we perform nonlinear analysis on the MDE.  We
first construct Poincar\'e maps of periodic solutions of the MDE to
examine asymptotic behavior of the moments and clarify that the
primary resonance encounters a saddle-node bifurcation.  We then
numerically continue the bifurcation point to obtain the saddle-node
bifurcation set that is in good agreement with the threshold of
occurrence of the slow convergence.

\section{Noise-induced Synchronization}

\subsection{Slow convergence of synchronization errors}

Let us consider the pair of identical piecewise linear
systems\cite{sss06,yoshida07}:
\begin{equation}
 \leqn{eq4-18}
  \begin{array}{l@{\;}l@{}}
   \ddot{x}+c\dot{x}+kG(x;\mu)&=Q+u(t),\\
   \ddot{y}+c\dot{y}+kG(y;\mu)&=Q+u(t).
  \end{array}
\end{equation}
where $Q$ is a preload, $u(t)$ is a random input, and the function
$G$ is a piecewise linear function defined by
\begin{math}
 G(x;\mu)= (x+\mu)-|x+\mu|+(x-\mu)+|x-\mu|
\end{math}
describing a linear spring with a dead zone of the width $2\mu$.  In
the previous work\cite{yoshida07}, we have examined the two types of
the input $u(t)$ such as the combination of the harmonic forcing and
white noise:
\begin{equation}
 \leqn{HR}
 u(t)=P\cos(\omega t)+sw(t),
\end{equation}
and the filtered noise: 
\begin{equation}
 \leqn{NB}
  \ddot{u}+2\zeta\omega_n\dot{u}+\omega_{n}^{2}u=sw(t),
\end{equation}
where $w(t)$ is the standard Gaussian white noise.  It is shown that
the synchronization system \reqn{eq4-18} produces the perfect
synchronization where the synchronization error vanishes
deterministically.

On the contrary, in the present work, we consider the random-phase
harmonic forcing as a purely random input in the following form:
\begin{equation}
 \leqn{RP}
 u(t)=P\cos\Big(\omega t + \rho B_t\Big)
\end{equation}
where $B_t$ is the standard Brownian motion.  We choose $c=0.04$,
$k=1.0$, $\mu=0.7$, $Q=0.3$, $P=0.2$, and $\rho=2\times10^{-5}$
which produce both the fast and slow convergence by sweeping
$\omega$ and can be implemented easily as a realistic mechanical
structure.

\Rfig{F1.SE} shows sample paths of synchronization errors for
$\omega=0.8$ (upper side) and $0.95$ (lower side) respectively.
There is significant difference in convergence time between the two
conditions. The convergence times are $t>20000$ and $t=118$
respectively.  This means that one spends time to converge more than
169 times as much as the other.

\begin{figure}[t]
\centering
  \includegraphics[width=.9\hsize]{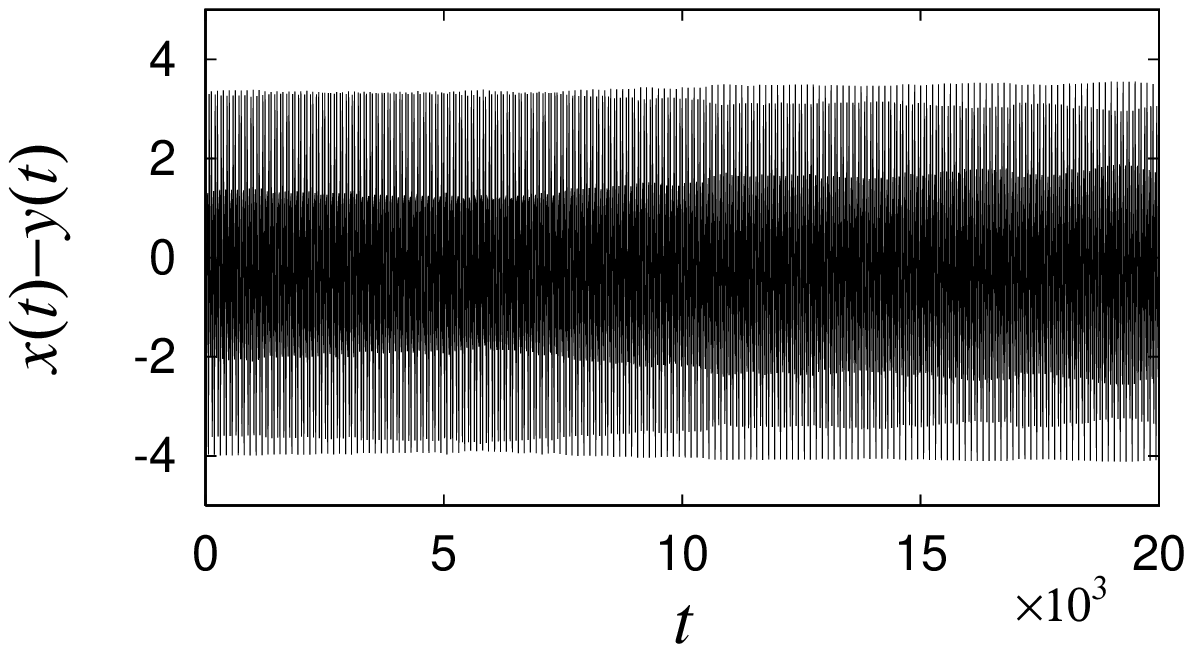}\\
  \includegraphics[width=.9\hsize]{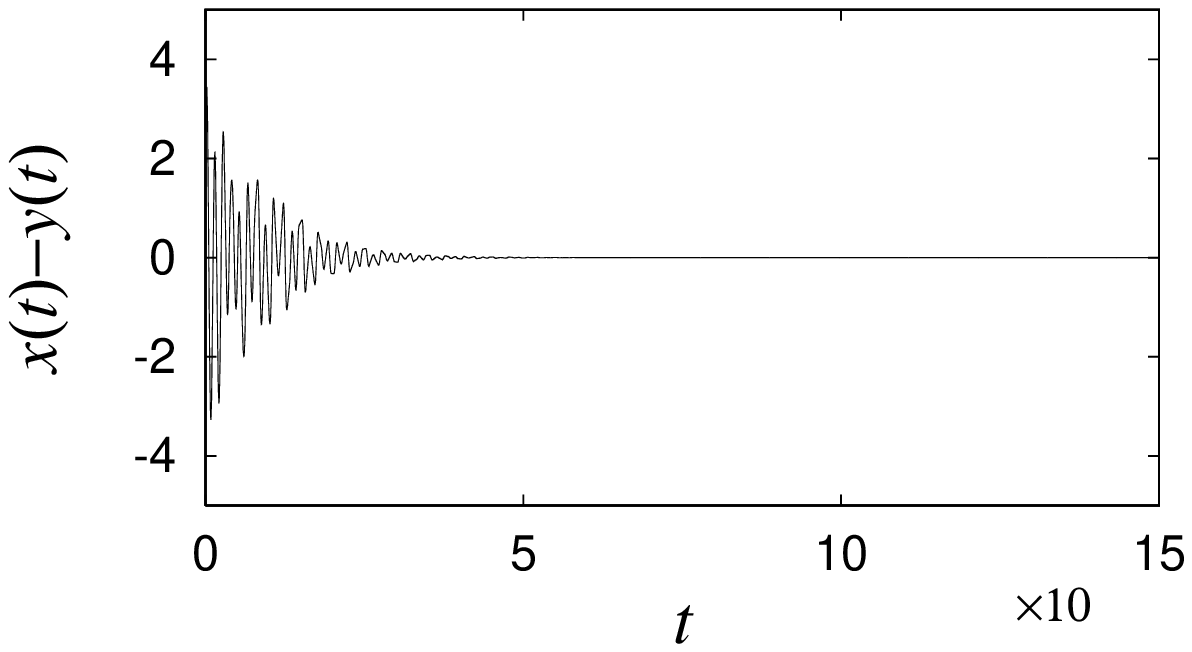}
  \caption{Sample paths of synchronization errors for $\omega=0.8$
  (upper) and $0.95$ (lower) respectively.}  
 \lfig{F1.SE}
\end{figure}

\subsection{Localization of probability densities}

One explanation of the slow convergence of the synchronization
errors is the transient localization of probability densities as
already discussed by the authors\cite{yoshida07}.

One of the pair \reqn{eq4-18}, say a single system, produces the
transient probability density $p(x_1,x_2,t)$ satisfying the following
Fokker-Plank-Kolmogorov (FPK) equation:
\begin{equation}
 \pd{p}{t} = \frac{\rho^2}{2}\pd{^2p}{x_3^2}-\pd{(\omega p)}{x_3}
 -\pd{(x_2p)}{x_1}
 -\pd{}{x_2}\Big(-cx_2-k\,G(x_1;\mu)+Q+P\cos x_3\Big)p
 \leqn{FPK:F1}
\end{equation}
where $x_1:=x$, $x_2:=\dot x$, $x_3:=B_t$.  To obtain direct solutions of these
equations, we employ Monte-Carlo techniques where the densities are
obtained from cumulative frequencies over the 10000 samples of
numerical solutions of one of the pair \reqn{eq4-18}.

\Rfig{hysteresis} shows dependency of the solution $p(x_1,x_2,t)$
of the FPK equation \reqn{FPK:F1} on the initial conditions
$p(x_1,x_2,0)$.  The small graphs represent the probability
densities at $t=60\times 2\pi/\omega$ starting from the initial
functions $p(x_1,x_2,0) = \delta(x_1-2)\delta(x_2-2)$ in the upper
side and $\delta(x_1+0.5)\delta(x_2)$ in the lower side
respectively.
It is obviously seen that the transient densities significantly
depend on the initial conditions within the region
$\omega\in[0.76,0.92]$.  In case of $\omega = 0.76$, both the initial
conditions (a) and (b) cause quite similar densities so that
dependency on the initial conditions is hardly found in this case.
As $\omega$ increases to $0.78$, the outer ring found at
$\omega=0.76$ is being replaced with the smaller ring for (a) while
the previous state maintains for (b).  Further increase of $\omega$
cause second change of the densities at $\omega=0.92$ where the
outer ring of (b) at $\omega=0.9$ is being replaced with the smaller
ring of (b) at $\omega=0.92$ while the densities of (a) maintains
between $\omega=0.9$ and $0.92$.

\begin{figure}[t]
\centering
 \vskip3ex \includegraphics[width=\hsize]{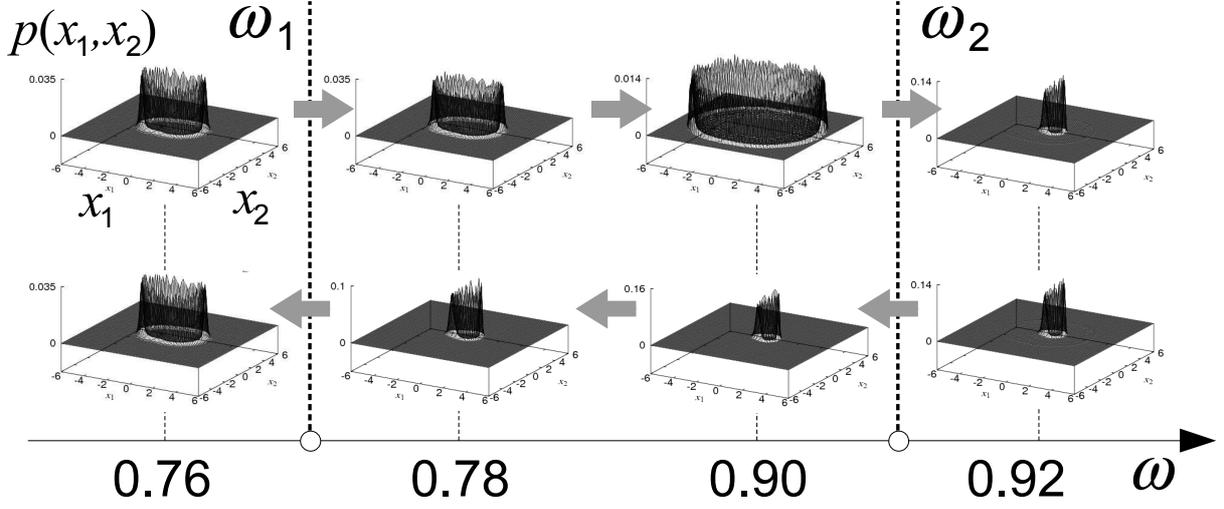}
 \caption{Transient isolation of probability densities at
 $t=60\times 2\pi/\omega$ comparable with deterministic hysteric
 jumps.}  \lfig{hysteresis}
\end{figure}

In the previous work\cite{yoshida07}, we have already discussed that
this dependency on initial conditions, in other words, the transient
isolation of the densities is one of the main reasons of the slow
convergence of the synchronization errors and can be detected by
hysteric jumps of statistical moments.  It will be shown that
similar results can be obtained even in this work considering the
random-phase forcing \reqn{RP}.

\subsection{Moment equations}

In order to evaluate the statistical moments of the synchronization,
we derive the moment differential equations (MDE) from
the FPK equation \reqn{FPK:F1} as follows.
\begin{align}
 &  \dot{m_{1}}=m_2,\;
   \dot{m_{2}}=-cm_{2}-k\ave{G} +Q +P\ave{\cos x_3},\notag\\
 &  \dot{m_{3}}=\omega,\;\dot{\sigma_{11}}=2\sigma_{12},\notag\\
 &  \dot{\sigma_{12}}=\sigma_{22}-c\sigma_{12}
       -k\ave{(x_1-m_1)G} + P\ave{x_1\cos x_3},\notag\\
 &  \dot{\sigma_{13}}=\sigma_{23}+\omega m_1,\notag\\
 &  \dot{\sigma_{22}}=-2c\sigma_{22}
       -2k\ave{(x_2-m_2)G} + 2P\ave{x_2\cos x_3},\notag\\
 &  \dot{\sigma_{23}}=-c\sigma_{23}
       -k\ave{(x_3-m_3)G} + P\ave{x_3\cos x_3} - \omega m_2,\notag\\
 &  \dot{\sigma_{33}}=2\omega m_3+\rho^2
 \leqn{mde}
\end{align}
where $m_i=\ave{x_i}$, $\sigma_{ij}=\ave{(x_i-m_i)(x_j-m_j)}$.  The
nonlinear averages are linearized by the standard statistically
equivalent techniques\cite{lin95,keijin84en}.  In our case, however,
$x_3=B_t$ is hardly assumed to be Gaussian because $B_t$ is
unbounded.  Therefore, utilizing the formula of the Brownian
motion\cite{Oksendal}:
\begin{equation}
 \ave{f\big(B_t\big)}=\frac{1}{\sqrt{2\pi t}}
  \int_{-\infty}^{\infty} f(x) e^{-x^2/2} dx,
\end{equation}
we perform the approximation
$\ave{\cos x_3} \approx e^{-\rho^2/2}\cos(m_3)$,
$\ave{x_3\cos x_3} \approx e^{-\rho^2/2}\sqrt{\omega m_3}\cos m_3 $.

In this setup, the moment equations \reqn{mde} can be regarded as
a periodically forced system (reducing $\dot m_3 = \omega$ to $m_3 =
\omega t$) where the third equation of \reqn{mde} is neglectable.
Furthermore, in order to calculate $\sigma_{11}$ and $\sigma_{22}$ only,
the sixth, eighth, and ninth equations are also neglectable because the
remaining equations are not coupled with them.  Therefore, we can
reduce the MDE \reqn{mde} into the following form:
\begin{align}
 &  \dot{m_{1}}=m_2,\notag\\
 &  \dot{m_{2}}=-cm_{2}-k\alpha_0(m_1,\sigma_{11})
       +Q +Pe^{\rho^2/2}\cos \omega t,\notag\\
 &  \dot{\sigma_{11}}=2\sigma_{12},\notag\\
 &  \dot{\sigma_{12}}=\sigma_{22}-c\sigma_{12}
       -k\alpha_1(m_1,\sigma_{11})\sigma_{11}
       + Pm_1 e^{\rho^2/2}\cos \omega t,\notag\\
 &  \dot{\sigma_{22}}=-2c\sigma_{22}
       -2k\alpha_1(m_1,\sigma_{11})\sigma_{12}
       + 2Pm_2 e^{\rho^2/2}\cos \omega t
 \leqn{mde:5}
\end{align}
where $\alpha_0(m_1,\sigma_{11})$ and $\alpha_1(m_1,\sigma_{11})$
are the statistically equivalent gains derived by
Sato\cite{keijin84en}:
\begin{align}
 \alpha_0 &= m_1 + \sqrt{\frac{\sigma_{11}}{2\pi}}
 \left( e^{-\frac{(m_1-\mu)^2}{2\sigma_{11}}}
       +e^{-\frac{(m_1+\mu)^2}{2\sigma_{11}}}\right)\notag\\
 +&
 \frac{m_1+\mu}{2}\erf\left(-\frac{m_1+\mu}{\sqrt{2\sigma_{11}}}\right)
 +\frac{m_1-\mu}{2}\erf\left(-\frac{m_1-\mu}{\sqrt{2\sigma_{11}}}\right),
 \\ \alpha_1 &= 1
 +\frac{1}{2}\erf\left(-\frac{m_1+\mu}{\sqrt{2\sigma_{11}}}\right)
 -\frac{1}{2}\erf\left(-\frac{m_1-\mu}{\sqrt{2\sigma_{11}}}\right)
\end{align}
where $\erf(\cdot)$ is the error function.  Therefore, in what
follows, we calculate $\sigma_{11}$ and $\sigma_{22}$ from
final expression \reqn{mde:5}.

\subsection{Jumps of moments}

\Rfig{SEL:F1} shows the peak-to-peak values of
$\ave{x_1^2}=\sigma_{22}$ calculated from the reduced MDE
\reqn{mde:5} and the mean convergence time\cite{yoshida07}:
\begin{align}
 \leqn{mct}
 \ave{T} := \frac{1}{K}\sum_{k=1}^{K} 
 \min\left\{
   n \left| \sup_{1\leq i,j\leq M} |X^i_k(n)-X^j_k(n)| < \epsilon\right.
 \right\}
\end{align}
where $n$ is a discrete time s.t. $t=n\Delta t$,
$\{X^m_k(n)\}_{m=1}^M$ is the $k$-th set of $M$ points composing
numerical probability densities, $k$ represents the $k$-th sample
path of the random excitation, and $\epsilon \ll 1$ is a criterion of
the convergence.  The moments are calculated as numerical solutions
of the moment equations \reqn{mde:5}. The solid and broken curves
indicate the forward and backward sweeps of the frequency $\omega$
respectively.  The mean convergence time $\ave{T}$ is estimated by
substituting numerical solutions of \reqn{eq4-18} starting from
$5\times 5$ uniform Cartesian grids on the region
$(x_1,x_2)\in[-10,10]\times[-10,10]$ into Eq.\reqn{mct} where
$M=5\times 5=25$, $K=100$, $\epsilon = 10^{-5}$.  The plots of
$\ave{T}$ are saturated at the given value $\ave{T}_{\max}=3000$.

\begin{figure}[t]
\centering
 \includegraphics[width=\hsize]{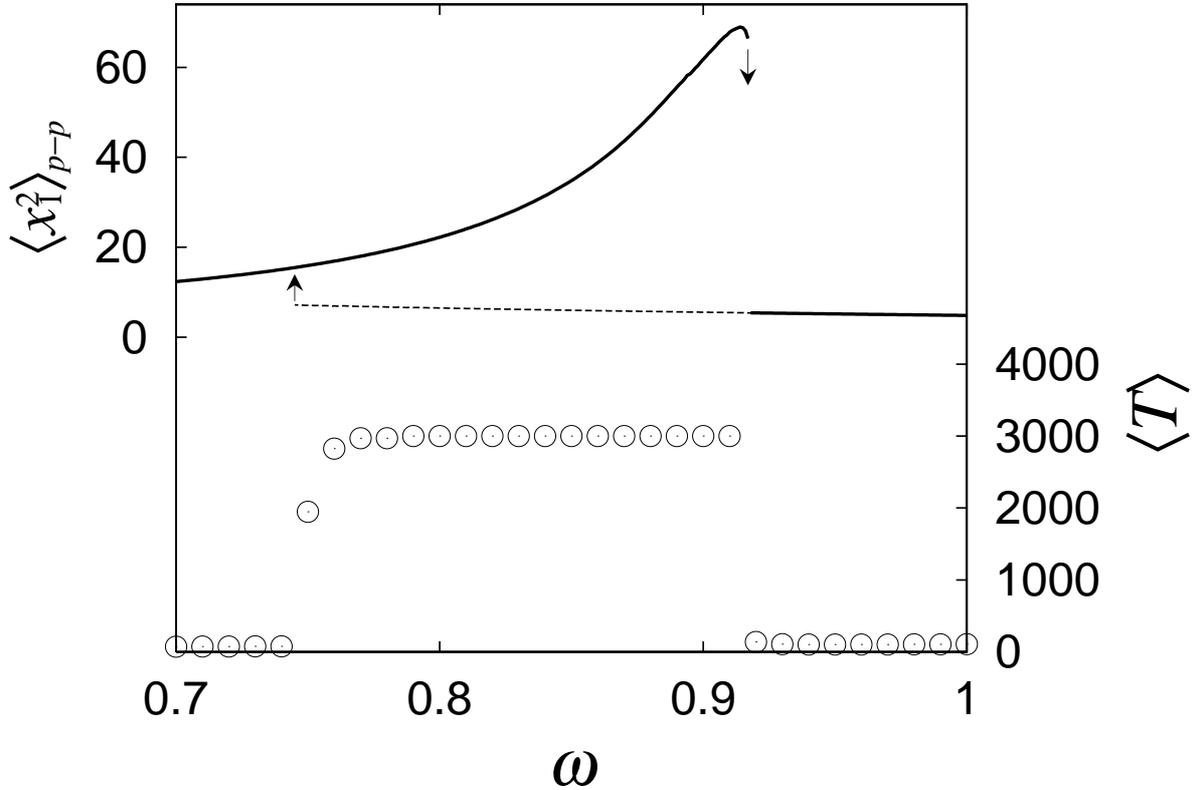}
 \caption{The peak-to-peak of the second moment $\langle
 x_1^2\rangle_{\text{p-p}}$ and the mean convergence time $\ave{T}$
 of the pair \reqn{eq4-18}.}  \lfig{SEL:F1}
\end{figure}

It is clearly seen from \rfig{SEL:F1} that the mean convergence time
$\ave{T}$ rapidly increase within the hysteric jumps of the second
moment of response $\ave{x_1^2}$.  This means that the hysteric
jumps of the moment in the statistical equivalent sense act as an
identifier to indicate the slow convergence of the synchronization
errors.  This means that the present result based on the random
phase forcing in Eq. \reqn{RP} agrees with the previous
result\cite{yoshida07} considering the other types of forcing 
in Eqs \reqn{HR} and \reqn{NB}.

\section{Bifurcation Analysis of Moments}

\subsection{Poincar\'e map of moments}
Let us now rewrite the MDE \reqn{mde:5} in the vector form:
\begin{equation}
 \dgv{q} = \gv{f}(\gv{q},t;p),\quad \gv{q}(0)=\gv{q}_0
 \leqn{mde:5v}
\end{equation}
where $t\in R$ is time, $p\in R$ is a free parameter, and $\gv
q=(m_1,m_2,\sigma_{11},\sigma_{12},\sigma_{22}):R\to R^5$,
$\gv{f}:R^5\times R\to R^5$. The solution of this initial value
problem can be written in the form:
\begin{equation}
 \gv{q}(t) = \varphi(t,\gv{q}_0)
\end{equation}
where $\varphi$ is the shift operator generated by the MDE. In our
case, the function $\gv{f}$ is supposed to be periodic in time,
\begin{equation}
 \gv{f}(\gv{q},t+\tau;p) = \gv{f}(\gv{q},t;p)
\end{equation}
where $\tau=2\pi/\omega$ is a period of the external forcing of the
MDE \reqn{mde:5v}.  Now, we provide a brief summary of how to
characterize periodic solutions of the MDE. See the
reference\cite{Holmes.pm} for details.  In order to reduce dimension,
the Poincar\'e map is defined by,
\begin{equation}
 T:R^5\to R^5:\gv q\mapsto \varphi(\tau,\gv q).
\end{equation}
In practice, a series of points generated by the Poincar\'e map is
obtained by $\gv q_{k}=T^k(\gv q_0)=\varphi(k\tau,\gv q_0)$~~
$(k=1,2,\cdots)$.  Plotting $\gv q_{k}$ is sometimes referred to as
the Poincar\'e plot.

Therefore, the $m$-periodic solutions passing through $\bar{\gv
q}\in R^5$ satisfies
\begin{equation}
 T^m(\bar{\gv q}) = \bar{\gv q},
  \leqn{poin}
\end{equation}
so that the periodic solutions can be regarded as the fixed points of
the mapping $T$.

\begin{figure}[t]
 \centering
 \includegraphics[width=\hsize]{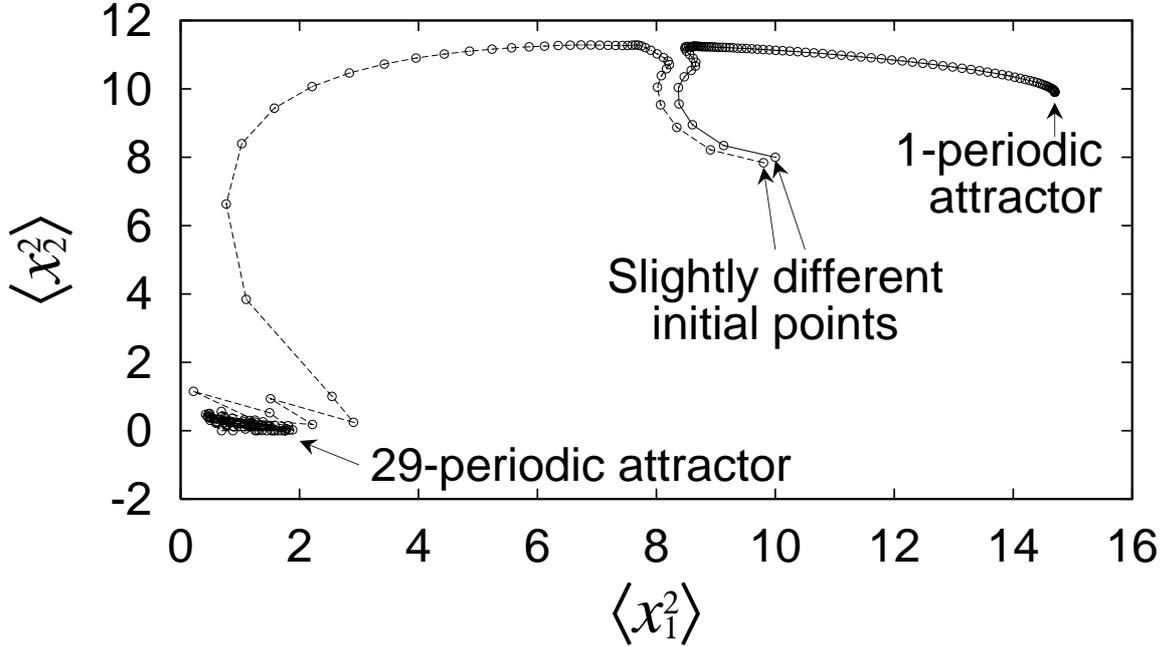}
 \caption{Poincar\'e plots of solutions of the MDE \reqn{mde:5}
 for $\omega=0.9$ starting from
 $\gv{q}(0)=\gv{a}=(0.5,2.5,10,2,8)$ and $0.98\gv{a}$.}
 \lfig{poin0.90}
\end{figure}

\Rfig{poin0.90} shows Poincar\'e plots of the solutions of the MDE
\reqn{mde:5} for $\omega=0.9$, which is inside of the hysteric
region of the response in \rfig{SEL:F1}, starting from the initial
points $\gv{q}(0)=\gv{a}=(0.5,2.5,10,2,8)$ and
$\gv{q}(0)=0.98\gv{a}$.  It is clearly seen that the initial point
$\gv{q}(0)=\gv{a}$ converges the 1-periodic attractor far from the
origin while $\gv{q}(0)=0.98\gv{a}$ converges to a 29-periodic
attractor near the origin. In this case, it seems that there is a
repeller near the initial points so that slight difference in
initial conditions results in distinct amplitudes of response as
already shown in \rfig{SEL:F1}.

On the other hand, \rfig{poin0.95} shows Poincar\'e plots for
$\omega=0.95$, which is outside of the hysteric region of the
response in \rfig{SEL:F1}, starting from the same initial conditions
as those in \rfig{poin0.90}.  It is obviously shown that the attractor far
from the origin vanishes and the periodic attractor near the origin
is replaced with a quasi-periodic attractor corresponding to a
solution with infinite period.  \Rfig{PDF:quasi} shows that
stationary probability density of the original equation
\reqn{eq4-18} under the condition corresponding to the
quasi-periodic attractor in \rfig{poin0.95} for $\omega=0.95$.
%\REV{
The two peaks of the density are caused by temporal switching
orbits of the sample paths fluctuating around bistable periodic orbits
of the original equation \reqn{eq4-18} under the deterministic limit.
%}
%This density implies that the sample paths undergo random switching in
%bistable manners.
%
It follows that the quasi-periodicity of the MDE does not reflect this
bistability of the sample paths.  One explanation of this disagreement
is the approximation errors of the nonlinear moments.

\begin{figure}[t]
 \centering \includegraphics[width=\hsize]{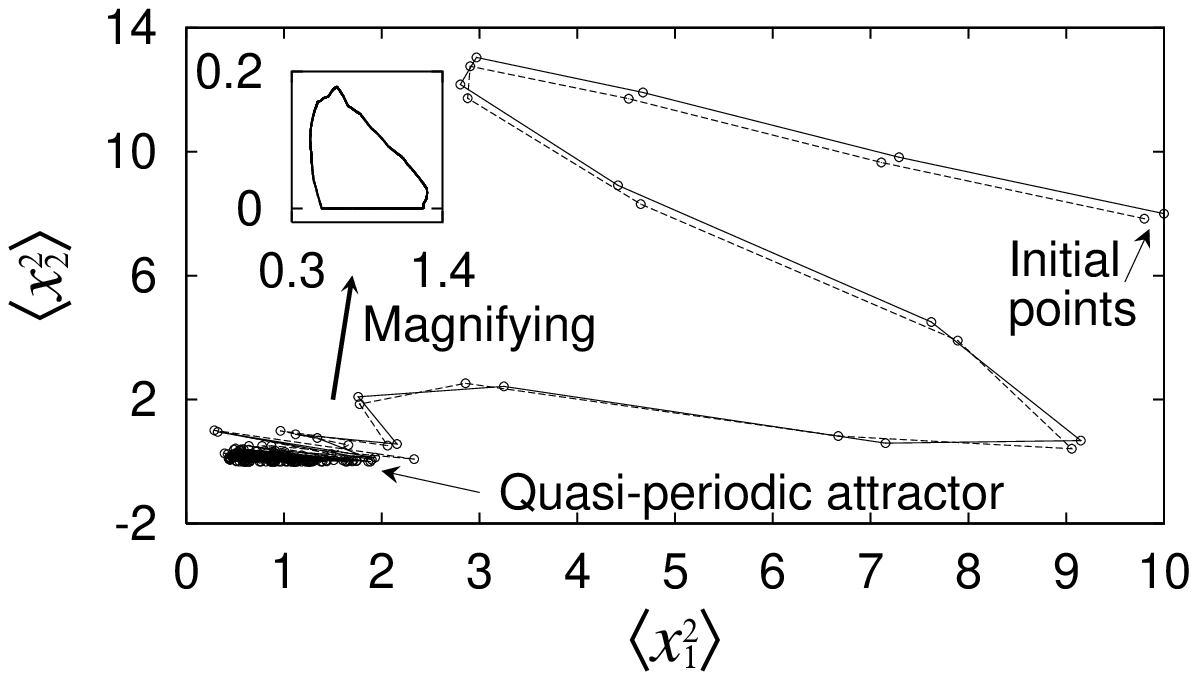}
 \caption{Poincar\'e plots of solutions of the MDE \reqn{mde:5}
 for $\omega=0.95$ starting from
 $\gv{q}(0)=\gv{a}=(0.5,2.5,10,2,8)$ and $0.98\gv{a}$.}
 \lfig{poin0.95}
 \medskip
%\end{figure}
%\begin{figure}[t]
 \centering \includegraphics[width=.8\hsize]{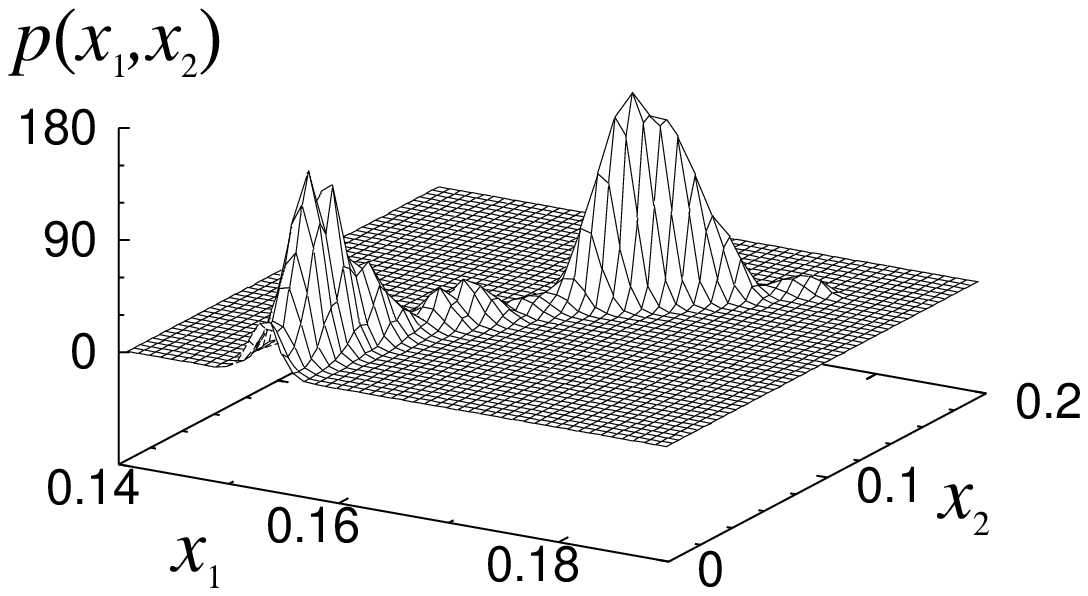}
 \caption{Stationary probability density of the original SDE
 \reqn{eq4-18} under the condition corresponding to the
 quasi-periodic attractor in \rfig{poin0.95} for
 $\omega=0.95$. Quasi-periodicity is hardly observed.}
 \lfig{PDF:quasi}
\end{figure}

From the above result, we can suppose that the primary
resonance (the large amplitude) in the response curve in
\rfig{SEL:F1} undergoes a saddle-node bifurcation at $\omega\approx
0.91$ because the pair of the repeller and the 1-periodic attractor
in \rfig{poin0.90} for $\omega=0.9$ vanishes in \rfig{poin0.95} for
$\omega=0.95$ while the off resonance undergoes other types of
bifurcation such as pitchfork bifurcations, Hopf bifurcations, and so
on.

In what follows, we focus on the primary resonance undergoing the
saddle-node bifurcation near $\omega=0.91$ to investigate how it
indicates the slow convergence of the synchronization errors.

\subsection{Numerical continuation of bifurcation points}

Let us summarize how to obtain bifurcation points numerically. See
Refs\cite{Holmes.pm,kawakami84} for more details.

It is mathematically proved that a matrix representation
$DT(\gv{q})$ of the linearization of the Poincar\'e map $T$ of the
ordinary differential equation \reqn{mde:5v} is given by,
\begin{equation}
 DT(\gv{q}):=J(\tau),\;\dot{J} = D\gv{f}(\gv{q},t;\lam) J,\; J(0) = I
  \leqn{DT}
\end{equation}
where $D\gv{f}(\gv{q},t;\lam)$ is the Jacobian matrix and $I$ is the
unit matrix. It can also be shown that the stabilities of the fixed
point $\bar{\gv q}$ of the Poincar\'e map $T$ is characterized by
the characteristic value $s$ satisfying the characteristic equation:
\begin{equation}
 P(s)=\left|sI-DT(\bar{\gv q})\right| = 0 \leqn{hyper}.
\end{equation}
From the basic theory of bifurcation\cite{kawakami84}, it is known that
saddle-node bifurcations occurs at $s=1$.

\begin{figure}[t]
 \centering \includegraphics[width=\hsize]{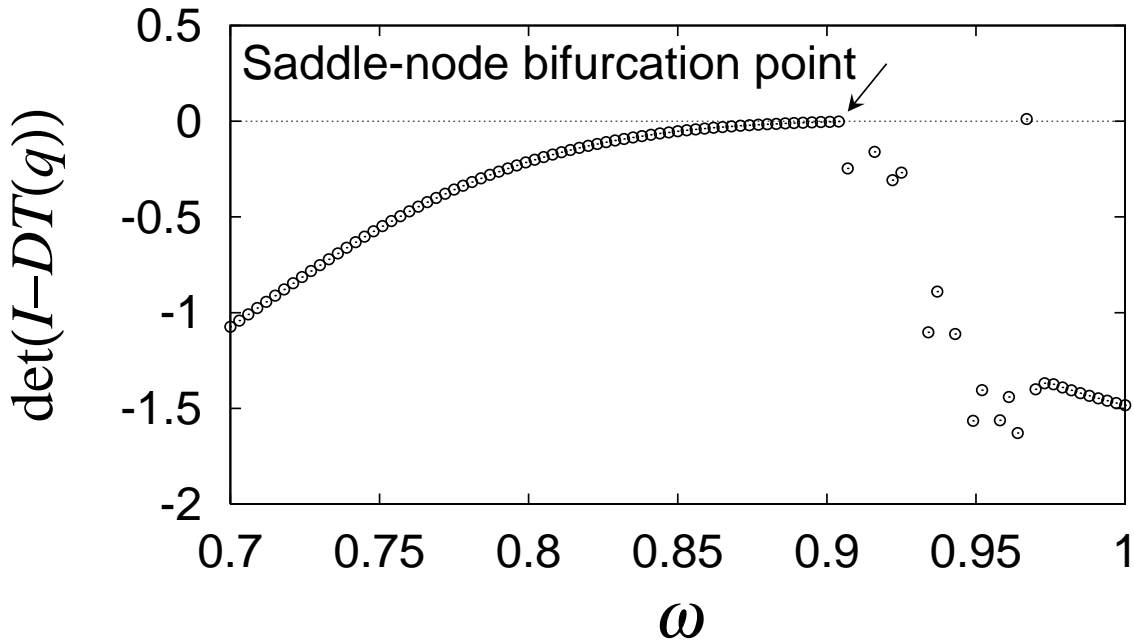}
 \caption{Characteristic polynomial $P(s)=|sI-DT(\gv{q})|$ for
 $s=1$ as a function of a free parameter $\omega$. $P(1)$ vanishes
 at the saddle-node bifurcation point.}
 \lfig{BIF:cp}
\end{figure}

\Rfig{BIF:cp} represents the value of characteristic polynomial for
$s=1$ calculated from Eqs \reqn{mde:5v} and \reqn{DT} for $c=0.04$.
It is clearly seen that the value $P(1)$ vanishes near $\omega=0.91$.
%\REV{
Since zeros of characteristic polynomials generally indicate
bifurcation points\cite{Holmes.pm},
%}
it is numerically proved that the primary resonance of the
second moment undergoes the saddle-node bifurcation near
$\omega=0.91$.

\subsection{Saddle-node bifurcation sets indicating the slow convergence}

In order to continue the bifurcation points in the parameter space
\cite{kawakami84}, we define new state vectors $\gv{p}:=(\bar{\gv
q},\lam) \in R^6$ and solve the equations \reqn{poin}
and \reqn{hyper} simultaneously:
\begin{equation}
 T^m(\bar{\gv q}) = \bar{\gv q},\quad
  P(s)=\left|sI-DT(\bar{\gv q})\right| = 0 \leqn{bifeqn}.
\end{equation}
The solution $\gv{p}$ is called a bifurcation set of the map
\reqn{poin}, and consequently yields a bifurcation set of the
corresponding periodic solution of the MDE \reqn{mde:5}.

\begin{figure}[t]
 \centering \includegraphics[width=\hsize]{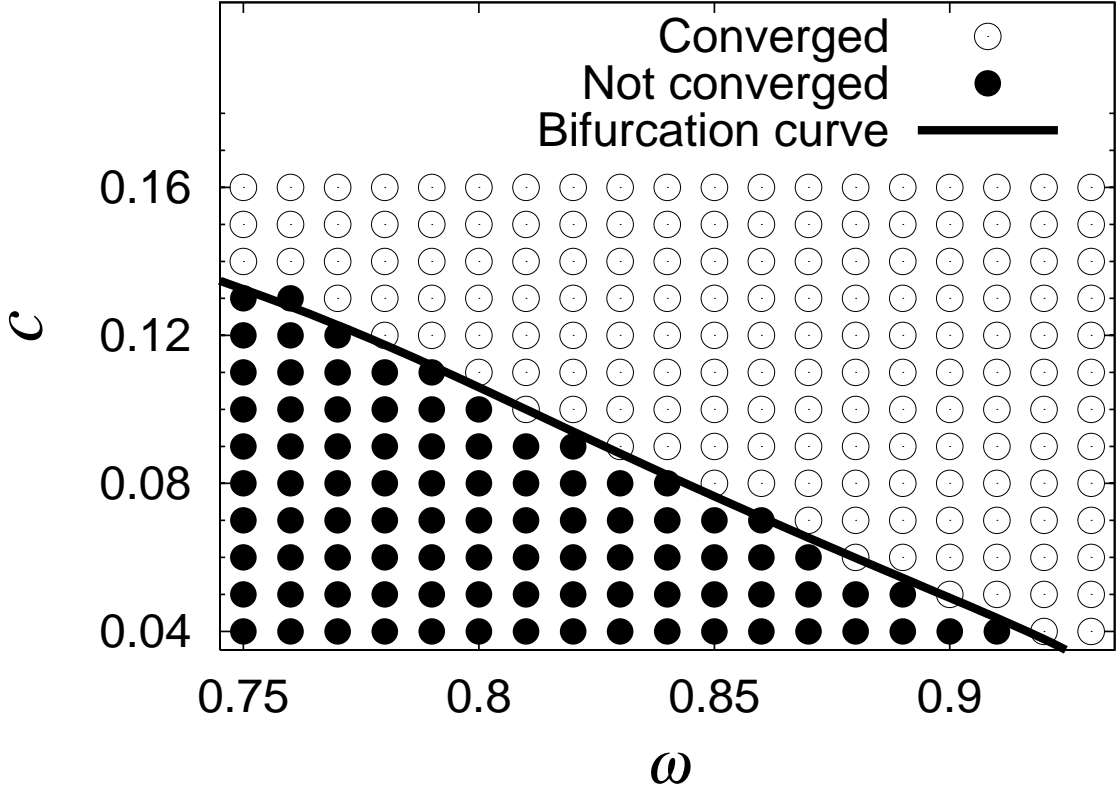}
 \caption{Numerical continuation of the saddle-node bifurcation
 point of periodic solutions of \reqn{mde:5} and mean convergence
 time $\ave{T}$ of the synchronization errors. The white circles are
 plotted for $\ave{T}<10^3$ and the black circles for otherwise
 respectively.}
 \lfig{BIF}
\end{figure}

\Rfig{BIF} shows the saddle-node bifurcation set in the
$(\omega,c)$-plane of periodic solutions of the MDE \reqn{mde:5} and mean
convergence time $\ave{T}$ of the synchronization errors.  In
practice, we consider the saddle-node bifurcation point in
\rfig{BIF:cp} for $c=0.04$ as the initial point and regard $\omega$
as the bifurcation parameter $\lam$. We then solve the
simultaneous equations \reqn{bifeqn} by standard Newton's methods.
Shifting the second parameter as $c=0.04+j\cdot\Delta c$
($j=1,2,\cdots$), we obtain a series of $\gv{p}$ on the
$(\omega,c)$-plane where $\Delta c$ is a small value chosen for a
good convergence of the Newton's methods.

The solid curve indicates the saddle-node bifurcation set of the MDE
\reqn{mde:5}.  The white and the black circles represent the fast
$\ave{T}<10^3$ and the slow $\ave{T}\geq 10^3$ convergences where
$\ave{T}$ is calculated from Monte-Carlo simulations of the
synchronization system \reqn{eq4-18}.
%\REV{
Along the bifurcation set, the solid curve, the characteristic
polynomial takes zero values at which the stability changes. In other
words, the changes of parameters crossing this curve transversely results
in the stability change of the system.
%}

It is very clear from the above result that the saddle-node
bifurcation set belonging to the primary resonance indicated by the
solid curve in \rfig{BIF} is in good agreement with the threshold of
occurrence of the slow convergence.

%%%%%%%%%%%%%%%%%%%%%%%%%%%%%%%%%%%%%%%%%%%%%%%%

\section{Conclusions}

We have investigated how to characterize the slow convergence of the
synchronization errors of the synchronization system which consists
of the pair of the piecewise linear systems subjected to the common
random excitation.  We first have demonstrated numerically that the
convergence speed of the synchronization errors significantly
depends on the parameter conditions.  It is shown that the slow
convergence is caused by the multimodal probability densities so
that it can be detected as the multi-valued solutions of the moment
differential equations.

We then have constructed Poincar\'e maps of periodic solutions of
the MDE to examine asymptotic behavior of the moments and clarify
that the primary resonance encounters saddle-node bifurcations.  We
have numerically continued the bifurcation point to obtain the
saddle-node bifurcation set of primary resonance of moments.  It is
clearly shown that the saddle-node bifurcation set belonging to the
primary resonance is in good agreement with the threshold of
occurrence of the slow convergence.

The above result leads to the conclusion that the saddle-node
bifurcation set of the moment equations provides an effective
identifier to detect the slow convergence of the synchronization errors.

%\bibliographystyle{sss}
%\bibliography{NIS}

%%%%%%%%%%%%%%%%%%%%%%%%%%%%%%%%%%%%%%%%%%%%%%%%%%%%%%%%%%%%%

\end{document}